\documentclass[11pt]{article}
\usepackage{mathrsfs}
\usepackage{amssymb}
\usepackage{amsmath}
\usepackage{amsfonts}
\usepackage{graphics}
\usepackage{epsfig}
\usepackage{enumerate}
\usepackage{color}

\newtheorem{theorem}{Theorem}[section]
\newtheorem{lemma}[theorem]{Lemma}
\newtheorem{definition}[theorem]{Definition}
\newtheorem{proposition}[theorem]{Proposition}
\newtheorem{corollary}[theorem]{Corollary}

\newenvironment{keywords}{{\bf Keywords: }}{}

\numberwithin{equation}{section}

\begin{document}

\author{Tianyang Nie$^{a,b,1}$, Aurel Rascanu$^{b,c,2}$ \\
%EndAName
{\small $^{a}$~School of Mathematics, Shandong University, Jinan,
Shandong 250100,
China}\\
{\small $^{b}$~Faculty of Mathematics, \textquotedblleft Alexandru
Ioan Cuza\textquotedblright\ University,}\\ {\small Bd. Carol I, no. 9-11,
Iasi 700506,
Romania}\\
{\small $^{c}$~\textquotedblleft Octav Mayer\textquotedblright\
Mathematics Institute of the Romanian Academy,}\\{\small Bd. Carol I, no.8,
Iasi 700506, Romania}}
\title{Direct and inverse images for fractional stochastic tangent sets and
applications }
\maketitle

\begin{abstract}
In this paper, we study direct and inverse images for fractional stochastic
tangent sets and we establish the deterministic necessary and sufficient
conditions that guarantee that the solution of a given stochastic
differential equation driven by the fractional Brownian motion evolves in
some particular sets $K$. As a consequence, a comparison theorem is obtain.
\end{abstract}

\noindent
\begin{keywords}
Stochastic Viability, Stochastic Differential Equations, Stochastic Tangent
Sets, Fractional Brownian Motion.
\end{keywords}

\let\oldthefootnote\thefootnote

\renewcommand{\thefootnote}

\footnote{%
E-mail addresses: tianyang.nie@uaic.ro, nietianyang@163.com; aurel.rascanu@uaic.ro.}

\let\thefootnote\oldthefootnote

\footnotetext[1]{%
Research partially supported by Marie Curie ITN Project, \textquotedblleft
Controlled Systems\textquotedblright , no.213841/2008.} \footnotetext[2]{%
Research partially supported by IDEAS project, no. 395/2007.}

\section{Introduction}
A general result on the existence and uniqueness of the solution for multidimensional, time dependent,
stochastic differential equations (SDEs) driven by a fractional Brownian motion (fBm) with Hurst parameter $H>\frac{1}{2}$
has been given by Nualart and R\u{a}\c{s}canu in~\cite{pl} using a techniques of the classical fractional calculus.

The notion of viable trajectories, used in the theory
of deterministic and stochastic differential equations, refers to those trajectories which
remain at any time in a fixed subset of the state space. The viability is to find necessary and sufficient conditions such that
a fixed subset is viable for the differential equation.  %Characterizing a subset and an equation in order for each solution
%to the equation starting from the subset to be viable in the subset is another
%problem, called the invariance problem.
In the theory of viable solutions the concept of the tangent sets and contingent sets play
a fundamental role. In fact, the pioneering theorem, proved in 1942 by
Nagumo, gives a criterion of the viability in terms
of contingent sets. Namely, the Nagumo theorem states that if $f$ is a
bounded, continuous map from a closed subset $K$ of $\mathbb{R}^m$ to $\mathbb{R}^m$, then a
necessary and sufficient condition such that $K$ is viable for the differential equation
\[x^{\prime}(t) = f(x(t)),  \quad x(0) = x_{0}\in K.
\]
is that
\[\langle f(x),p\rangle\leq0,\quad\forall\text{~}x\in K\text{ and }%
\forall\text{~}p~\;\text{a normal vector at}\;K\text{ in }x.
\]

Various generalizations of the Nagumo theorem provide viability conditions
in terms of contingent cones (see for instance ~\cite{pa} Th. 1, p. 191).
Viability and invariance with respect to It\^{o} equations have been investigated
first by J.-P. Aubin and G. Da Prato in~\cite{pc}.
%The stochastic contingent set defined in that paper is an adaptation of the deterministic concept.
Criterions for the viability and invariance of closed and convex subset of $\mathbb{R}^m$, given
in~\cite{pc}, are expressed in terms of stochastic contingent sets. Their results were generalized to arbitrary subsets
(which can also be time-dependent and random) in~\cite{pj}.

Another approach has been developed by Buckdahn, Peng, Quincampoix, Rainer and R\u{a}\c{s}canu in~\cite{pe},~\cite{pf},~\cite{pg},~\cite{ph}.
The main point of their work consist in proving that the viability property for SDE and
also for backward SDE holds true if and only if the square of the distance to the constraint sets is a viscosity supersolution(subsolution)
of the Hamilton-Jacobi-Bellman equation associated.

With respect to the SDE driven by fBm,  I. Ciotir and {A.R\u{a}\c{s}canu} proved a type of Nagumo Theorem on viability properties of close bounded
subsets with respect to a stochastic differential equation driven by fractional Brownian motion in~\cite{pi}.

Conditions expressed by stochastic contingent sets which are given in~\cite{pi} are general but unfortunately
not easy to check and the aim of the present paper is to give
checkable conditions for general stochastic differential equation driven by the fractional Brownian
motion and some particular sets $K$.

Studying from~\cite{pi}, we find the deterministic necessary and sufficient conditions that guarantee
that the solution of a  stochastic differential equation driven by the fractional Brownian motion $B^{H}$ with Hurst parameter $\frac{1}{2}<H<1$
(in short: f-SDE), $\mathbb{P}$-$a.s. \omega\in\Omega$
\[
X_{s}^{t,x}=x+{\displaystyle\int_{t}^{s}}b(r,X_{r}^{t,x})dr+{\displaystyle\int_{t}^{s}}\sigma(r,X_{r}^{t,x})dB_{r}^{H}%
,\quad s\in[t,T],
\]
$(t,x)\in[0,T]\times\mathbb{R}^d$ envolves in some  particular sets $K$ i.e. under which it holds that for all $t\in[0,T]$
and for all $x\in K$:
\[
X_{s}^{t,x}\in K  \quad a.s. \omega\in\Omega, \quad\forall s\in[t,T]  .
\]
Here

\begin{itemize}
\item $B^{H}=\left\{  B^{H}_{t},t\geq0\right\}  $ is a
fractional Brownian motion with Hurst parameter $\frac{1}{2}<H<1$, and the integral with respect to $B^{H}$ is a pathwise Riemann-Stieltjes integral;

\item
$b(t,x) :\left[  0,T\right]  \times
\mathbb{R}^{d}\mathbb{\rightarrow}\mathbb{R}^{d}$ and
$\sigma(t,x) :\left[  0,T\right]  \times\mathbb{R}%
^{d}\mathbb{\rightarrow R}^{d}$ are continuous functions.
\end{itemize}

The characterization of viability of $K$ is obtained through the study of the direct and inverse images for fractional stochastic tangent sets.
This idea comes from~\cite{pc}. In fact we extend the direct and inverse images of stochastic tangent sets to the fractional form and
using our main theorem \ref{th3.2},
we character the viability of some particular sets $K$ with the conditions on $b$ and $\sigma$ and we also obtain a comparison theorem.

We now explain how the paper is organized. In the second section,
we recall some classical definitions and the assumptions on the coefficients supposed to hold.
we also recall the main result in~\cite{pi}, which we will use later.
In section 3 we state our main result and some applications are given.
The section 4 is devoted to the proof of the main result and section 5 is for the proof of a general comparison theorem.

\section{Preliminaries}
Consider the equation on $\mathbb{R}^{d}$
\begin{equation}\label{equ2.2}
X_{s}=X_{0}+{\displaystyle\int_{0}^{s}}b(r,X_{r})dr+{\displaystyle\int_{0}^{s}}\sigma(r,X_{r})dB_{r}^{H}
,\quad s\in[0,T],
\end{equation}
\begin{itemize}
\item $B^{H}=\left\{  B^{H}_{t},t\geq0\right\}  $ is a
fractional Brownian motion defined on a complete probability space $(\Omega
,\mathcal{F},\mathbb{P)};$ with Hurst parameter $\frac{1}%
{2}<H<1$, and the integral with respect to $B^{H}$ is a pathwise Riemann-Stieltjes integral;

\item $X_{0}$ is a $d$ - dimensional random variable.

\item $b  : [0,T]  \times
\mathbb{R}^{d}\mathbb{\rightarrow}\mathbb{R}^{d}$, $\sigma  : [0,T]  \times\mathbb{R}%
^{d}\mathbb{\rightarrow R}^{d}$ are continuous functions.
\end{itemize}

Remark that the fractional Brownian motion has the following property:

For every $0<\varepsilon<H$ and $T>0$ there exists a positive random variable $\eta_{\varepsilon,T}$ such that
$\mathbb{E}(|\eta_{\varepsilon,T}|^{p})<\infty$, for all $p\in[1,\infty)$ and for all $s,t\in[0,T]$
\[
|B^{H}(t)-B^{H}(s)|\leq\eta_{\varepsilon,T}|t-s|^{H-\varepsilon}      \quad a.s.
\]
And from~\cite{pm} proposition 1.7.1(see also in \cite{pn}), we have for every $t_{0}\in[0,+\infty)$,
\begin{equation}\label{equ2}
\mathbb{P}\left\{\limsup_{t\rightarrow t_{0},~t\ge
t_{0}}\Big|\frac{B^{H}(t)-B^{H}(t_{0})}{t-t_{0}}\Big|=+\infty\right\}=1
\end{equation}
Using the same method we can easily proof that
\begin{equation}\label{equ3}
\mathbb{P}\left\{\limsup_{t\rightarrow t_{0},~t\ge
t_{0}}\frac{B^{H}(t)-B^{H}(t_{0})}{t-t_{0}}=+\infty\right\}=\mathbb{P}\left\{\liminf_{t\rightarrow
t_{0},~t\ge
t_{0}}\frac{B^{H}(t)-B^{H}(t_{0})}{t-t_{0}}=-\infty\right\}=\frac{1}{2}
\quad .
\end{equation}
\subsection{Assumptions and Notations}

For the function and coefficients appearing in the equation (\ref{equ2.2}),
we make the following standard assumptions which we will use throughout the paper:
\begin{description}
\item[$\left(  \mathbf{H}_{1}\right)  $]$\sigma(t,x)$ is differentiable in $x\in\mathbb{R}^{d}$,
and there exist some constants $\beta, \delta, 0<\beta,\delta\leq1$, and for every $R>0$ there exists $M_{R}>0$
such that the following properties hold for all $t\in[0,T]$,
\[
(H_{\sigma}):\;\left\{
\begin{array}
[c]{rl}%
i)\quad & \left\vert \sigma(t,x)-\sigma(s,y)\right\vert \leq M_{0}\left(
|t-s|^{\beta}+|x-y|\right)  ,\quad\forall x,y\in\mathbb{R}^{d},\medskip\\
ii)\quad & |\nabla_{x}\sigma(t,y)-\nabla_{x}\sigma(s,z)|\leq M_{R}\left(
|t-s|^{\beta}+|y-z|^{\delta}\right)  ,\quad\forall\left\vert y\right\vert
,\left\vert z\right\vert \leq R,
\end{array}
\right.
\]
where $\nabla_{x}\sigma(t,x)=(\nabla_{x}\sigma^{i}(t,x))_{i=\overline{1,d}}$ and
\[
|\nabla_{x}\sigma(t,x)|^{2}=\sum_{l=1}^{d}\sum_{i=1}^{d}|\partial_{x_{l}}\sigma^{i}(t,x)|^2
\]
Remark that for all $x\in\mathbb{R}^{d}$
\[
|\sigma(t,x)|\leq|\sigma(0,0)|+M_{0}(|t|^{\beta}+|x|)\leq M_{0,T}(1+|x|)
\]
where $M_{0,T}=|\sigma(0,0)|+M_{0}+M_{0}T$.

Let
\[
\alpha_{0}=\min\left\{  \frac{1}{2},\beta,\frac{\delta}{1+\delta}\right\}  .
\]
\item[$\left(  \mathbf{H}_{2}\right)  $] There exist $\mu\in(1-\alpha_{0},1]$
and for every $R\geq0$ there exists $L_{R}>0$ such that the following properties hold for all $t\in\left[
0,T\right]  ,$
\[
(H_{b}):\;\left\{
\begin{array}[c]{rl}
i)\quad & \left\vert b(r,x)-b(s,y)\right\vert \leq L_{R}\left(|r-s|^{\mu}+|x-y|\right)  ,\quad\quad\forall\left\vert x\right\vert
,\left\vert y\right\vert \leq R,\medskip\\
ii)\quad & \left\vert b(t,x)\right\vert \,\leq\,L_{0}(1+|x|),\quad\forall
x\in\mathbb{R}^{d}.
\end{array}
\right.
\]
\end{description}
Finally, we introduce some notations which will be used later.

Let $d,k\in\mathbb{N}^{*}$. Given a matrix $A=(a^{i,j})_{d\times k}$
and a vector $y=(y^{i})_{d\times 1}$, we denote
$|A|^{2}=\sum_{i,j}|a^{i,j}|^{2}$ and $|y|=\sum_{i}|y^{i}|^{2}$.
\medskip

Let $t\in\left[  0,T\right]  $ be fixed. Denote

\begin{itemize}
\item $W^{\alpha,\infty}(t,T;\mathbb{R}^{d}),$ $0<\alpha<1,$ the space of
continuous functions $f:[t,T]\rightarrow\mathbb{R}^{d}$ such that%
\[
\left\Vert f\right\Vert _{\alpha,\infty;\left[  t,T\right]  }:=\sup
_{s\in\lbrack t,T]}\left(  |f(s)|+\int_{t}^{s}\dfrac{\left\vert f\left(
s\right)  -f\left(  r\right)  \right\vert }{\left(  s-r\right)  ^{\alpha+1}%
}dr\right)  <\infty.
\]
An equivalent norm can be defined by%
\[
\left\Vert f\right\Vert
_{\alpha,\lambda;\left[  t,T\right]  }:=\sup_{s\in\lbrack
t,T]}e^{-\lambda s}\left(  |f(s)|+%
%TCIMACRO{\dint _{t}^{s}}%
%BeginExpansion
{\displaystyle\int_{t}^{s}}
%EndExpansion
\dfrac{\left\vert f\left(  s\right)  -f\left(  r\right)  \right\vert }{\left(
s-r\right)  ^{\alpha+1}}dr\right)\quad\forall\lambda\geq0.
\]

\item $\tilde{W}^{1-\alpha,\infty}(t,T;\mathbb{R}^{d}),$ $0<\alpha<\frac{1}{2}.$ the space of
continuous functions $g:[t,T]\rightarrow\mathbb{R}^{d}$ such that

\[
\left\Vert g\right\Vert _{\tilde{W}^{1-\alpha,\infty}(t,T;\mathbb{R}^{d}%
)}:=\left\vert g\left(  t\right)  \right\vert +\sup_{t<r<s<T}\left(
\frac{|g(s)-g(r)|}{(s-r)^{1-\alpha}}+\int_{r}^{s}\frac{|g(y)-g(r)|}%
{(y-r)^{2-\alpha}}dy\right)  <\infty.
\]

\item $C^\mu([t,T];\mathbb{R}^d)$, $0<\mu<1$, the space of $\mu$-H\"{o}lder continuous functions $f:[t,T]\to\mathbb{R}^d$, equipped with the norm
\[\left\Vert f\right\Vert
_{\mu;\left[  t,T\right]  }:=\left\Vert  f\right\Vert
_{\infty;\left[  t,T\right]  }+\sup_{t\leq r<s\leq T}\dfrac{\left\vert f\left(  s\right)  -f\left(  r\right)  \right\vert }{\left(
s-r\right)  ^{\mu}}<\infty
\]
where $\left\Vert  f\right\Vert_{\infty;\left[  t,T\right]}:=\sup_{s\in[t,T]}\lvert f(s)\rvert$. We have, for all $0<\varepsilon<\alpha$
\[C^{\alpha+\varepsilon}([t,T];\mathbb{R}^d)\subset W^{\alpha,\infty}(t,T;\mathbb{R}^{d})
\]

\end{itemize}

\begin{itemize}
\item $W^{\alpha,1}(t,T;\mathbb{R}^{d})$ the space of measurable functions $f$
on $[t,T]$ such that

\[
\left\Vert f\right\Vert _{\alpha,1;\left[  t,T\right]  }:=
\int_{t}^{T}\left[  \frac{|f(s)|}{\left(  s-t\right)  ^{\alpha}}+\int_{t}%
^{s}\frac{|f(s)-f(y)|}{(s-y)^{\alpha+1}}dy\right]  ds<\infty.
\]
Clearly
\[
W^{\alpha,\infty}(t,T;\mathbb{R}^{d})\subset W^{\alpha,1}(t,T;\mathbb{R}%
^{d}).
\]
\end{itemize}
\subsection{Generalized Stieltjes integral}
Denoting
\[
\Lambda_{\alpha}(g;\left[  t,T\right]  ):=\frac{1}%
{\Gamma(1-\alpha)}\sup_{t<r<s<T}\left\vert \left(  D_{s-}^{1-\alpha}%
g_{s-}\right)  (r)\right\vert .
\]
where
\[
\Gamma(\alpha)=\int_{0}^{\infty}s^{\alpha-1}e^{-s}ds
\]
is the Gamma function and
\[(D_{s-}^{1-\alpha}g_{s-})(r)=\frac{e^{i\pi(1-\alpha)}}{\Gamma(\alpha)}\left(
\frac{g(s)-g(r)}{(s-r)^{1-\alpha}}+(1-\alpha)\int_{r}^{s}\frac{g(r)-g(y)}%
{(y-r)^{2-\alpha}}dy\right)1_{(t,s)}(r).
\]
we have
\[
\Lambda_{\alpha}(g;\left[  t,T\right]  )\leq\frac{1}
{\Gamma(1-\alpha)\Gamma(\alpha)}\left\Vert g\right\Vert _{\tilde{W}^{1-\alpha,\infty}(t,T;\mathbb{R}^{d}
)}
\]
Note that
\[
\Lambda_{\alpha}(g;\left[  t,T\right]  )\leq\Lambda_{\alpha}(g;\left[  0,T\right]  )\Big(:=\Lambda_{\alpha}(g)\Big).
\]
We also introduce the notation
\[
(D_{t+}^{\alpha}f)(r)=\frac{1}{\Gamma(1-\alpha)}\left(
\frac{f(r)}{(r-t)^{\alpha}}+\alpha\int_{t}^{r}\frac{f(r)-f(y)}%
{(r-y)^{\alpha+1}}dy\right)1_{(t,T)}(r).
\]
\begin{definition}\label{def2.1}
Let $0<\alpha<\frac{1}{2}$. If $f\in W^{\alpha,1}(t,T;\mathbb{R}^{d\times k})$ and
$g\in\tilde{W}^{1-\alpha,\infty}(t,T;\mathbb{R}^{k})$, then defining
\[
\int_{t}^{s}f\left(  r\right)  dg\left(  r\right)  :=\left(
-1\right)  ^{\alpha}\int_{t}^{s}\left(  D_{t+}^{\alpha}f\right)  (r)\left(
D_{s-}^{1-\alpha}g_{s-}\right)  \left(  r\right)  dr.
\]
the integral $%TCIMACRO{\dint _{t}^{T}}%
%BeginExpansion
{\displaystyle\int_{t}^{s}}
%EndExpansion
fdg
$ exists for all $s\in[t,T]$ and
\begin{eqnarray*}
\left\vert
%TCIMACRO{\dint _{t}^{T}}%
%BeginExpansion
{\displaystyle\int_{t}^{T}}
%EndExpansion
f\left(  r\right)  dg\left(  r\right)  \right\vert
& \leq &
\sup_{t\leq r<s\leq T}\lvert(D_{s-}^{1-\alpha}g_{s-})(r)\rvert
%TCIMACRO{\dint _{t}^{T}}%
%BeginExpansion
{\displaystyle\int_{t}^{T}}
%EndExpansion
\lvert(D_{t+}^{\alpha}f)(s)ds\rvert{}
\nonumber\bigskip \\
& \leq & {}
\Lambda_{\alpha
}(g;\left[  t,T\right]  )\left\Vert f\right\Vert _{\alpha,1;\left[
t,T\right]  .}
\end{eqnarray*}
\end{definition}
It is known that when $H\in(\frac{1}{2},1)$ and $1-H<\alpha<\frac{1}{2}$, then the random variable
\[
G=\Lambda_{\alpha}(
B^{H})=\frac{1}{\Gamma(1-\alpha)}\sup_{t<s<r<T}|(D_{r-}^{1-\alpha}B_{r-})(s)|
\]
has moments of all order. As a consequence, if $u=\{u_{t}, t\in[0,T]\}$ is a
stochastic process whose trajectories belong to the space $W^{\alpha,1}(t,T;\mathbb{R}^{d})$, with $1-H<\alpha<\frac{1}{2}$,
the pathwise integral $
%TCIMACRO{\dint _{t}^{T}}%
%Beginexpansion
{\displaystyle\int_{0}^{T}}
%Endexpansion
u_{s}dB_{s}^{H}$ exists in the sense of Definition \ref{def2.1} and we have the estimate
\[
\Big|
%TCIMACRO{\dint _{t}^{T}}%
%BeginExpansion
{\displaystyle\int_{0}^{T}}
%EndExpansion
u_{s}dB^{H}_{s}\Big|
\leq G\|u\|_{\alpha,1}.
\]
This is the reason why in the SDE (\ref{equ2.2}) the integral with respect to $B^{H}$ is a pathwise Riemann-Stieltjes integral.

D. Nualart and\ A. R\u{a}\c{s}canu have proved in~\cite{pl}
that under the assumpations $\left(
\mathbf{H}_{1}\right)  $ and $\left(  \mathbf{H}_{2}\right)  ,$ with $\beta
>1-H$ and $\delta>\frac{1}{H}-1$ the SDE
\begin{equation*}
X_{s}^{t,\xi}=\xi+\int_{t}^{s}b(r,X_{r}^{t,\xi})dr+\int_{t}^{s}\sigma\left(
r,X_{r}^{t,\xi}\right)  dB^{H}_{r},\,\;s\in\left[  t,T\right]  ,
\end{equation*}
 has a unique solution $X^{t,\xi}\in L^{0}\left(  \Omega
,\mathcal{F},\mathbb{P\,};W^{\alpha,\infty}(t,T;\mathbb{R}^{d})\right)
,$ for all $\alpha\in\left(  1-H,\alpha_{0}\right)  .$
Moreover, for $\mathbb{P}$-almost all $\omega\in\Omega$,
$X\left(  \omega,\cdot\right) \in C^{1-\alpha}\left(  0,T;\mathbb{R}^{d}\right)  .$

\subsection{Fractional Viability}

In this subsection we recall the notion of the viability property for SDE driven by fractional Brownian motion.
On the other hand we will present the main result of~\cite{pi} which is very useful for our results.

Consider the stochastic differential equation driven by fractional Brownian motion $B^{H}$ with Hurst parameter $\frac{1}{2}<H<1$,
\begin{equation}\label{equ2.3}
X_{s}^{t,x}=x+{\displaystyle\int_{t}^{s}}b(r,X_{r}^{t,x})dr+{\displaystyle\int_{t}^{s}}\sigma(r,X_{r}^{t,x})dB_{r}^{H}
,\quad s\in[t,T].
\end{equation}

\begin{definition}
Let $\mathcal{K}=\{  K(t)  :t\in[0,T]\}$
be a family of subsets of  $\mathbb{R}^{d}$. We say that $\mathcal{K}$ is viable (weak invariant) for the equation
(\ref{equ2.3}) if, starting at any time $t\in[0,T]$ and
from any point $x\in K(t)$, there exists at least one of its solutions $\{
X_{s}^{t,x}:s\in[t,T]\}$ which satisfies
\[X_{s}^{t,x}\in K( s) \quad \text{for all}\quad  s\in[t,T].
\]

\end{definition}
\begin{definition}
The family $\mathcal{K}$ is invariant (strong invariant) for the equation
(\ref{equ2.3}) if, for any $t\in[0,T]  $ and for any starting point $x\in K(t)$, all solutions
$\{  X_{s}^{t,x}:s\in[t,T]  \}  $ of the fractional stochastic differential equation
(\ref{equ2.3}) have the property
\[
X_{s}^{t,x}\in K(s)  \quad \text{for all} \quad s\in[t,T]  .
\]
\end{definition}
Remark that, in the case when the equation has a unique solution (which is the case for equation (\ref{equ2.3})
under the assumptions $(\mathbf{H}_{1})$ and $(\mathbf{H}_{2})$), viability is equivalent to invariance.
\smallskip

Assuming that the mappings $b$ and $\sigma$ from the equation (\ref{equ2.3}) satisfy $(\mathbf{H}_{1})$ and  $(\mathbf{H}_{2})$.
\begin{definition}
Let $t\in\left[  0,T\right]  $ and
$x\in K\left(  t\right)  .$ Let $\frac{1}{2}<1-\alpha<H.$
\medskip

The $\left(  1-\alpha\right)  $-fractional $B^{H}$-contingent set to
$K\left(  t\right)  $ in $x$ is the set of the pairs $(u,v)$, such
that there exist random variable $\bar{h}=\bar{h}^{t,x}>0$ and a
stochastic process $Q=Q^{t,x}:\Omega\times\left[  t,t+\bar{h}\right]
\rightarrow \mathbb{R}^{d}$, and for every $R>0$ with $\left\vert
x\right\vert \leq R$ there exist two random variables
$H_{R},\tilde{H}_{R}>0$ independent of $(t,\bar{h})$ and a constant
$\gamma=\gamma_{R}(\alpha,\beta)\in(0,1)$ such that for all
$s,\tau\in[t,t+\bar{h}]$, $\mathbb{P}$-a.s.
\[
\left\vert Q\left(  s\right)  -Q\left(  \tau\right)  \right\vert \leq
H_{R}\left\vert s-\tau\right\vert ^{1-\alpha},\quad\left\vert
Q\left(  s\right)  \right\vert \leq\tilde{H}_{R}\left\vert s-t\right\vert
^{1+\gamma}%
\]
and
\[
x+\left(  s-t\right) u  +v  \left[
B_{s}^{H} -B_{t}^{H} \right]  +Q\left(  s\right)  \in
K\left(  s\right)  ,
\]
where the constants $H_{R}$, $\tilde{H}_{R}%
$ depend only on $R$, $L_{R}%
$, $M_{0,T}$,$M_{0}$, $L_{0}$,
$T$, $\alpha$, $\beta$, $\Lambda_{\alpha}\left(
B^{H}\right)  $.
\end{definition}

\begin{definition}\label{def2.5}
Let $t\in\left[  0,T\right]  $ and
$x\in K(t).$ Let $\frac{1}{2}<1-\alpha<H.$

The $\left(  1-\alpha\right)  $-fractional $B^{H}$-tangent set
 to $K(t)  $ in $x$, denoted by $S_{K(t)}(t,x)$, is the set of the pairs $(u,v)$, such that there exist random variable $\bar{h}=\bar{h}^{t,x}>0$
and two stochastic process
\[
\begin{array}
[c]{ll}
U=U^{t,x}:& \left[  t,t+\bar{h}\right]  \rightarrow\mathbb{R}^{d} ,~\ U(t)=0
\\
V=V^{t,x}:& \left[  t,t+\bar{h}\right]  \rightarrow\mathbb{R}^{d} ,~\ V(t)=0
\end{array}
\]
and for every $R>0$ with $\left\vert x\right\vert \leq R$ there
exsit two random variables $D_{R},\tilde{D}_{R}>0$ independent of
$(t,\bar{h})$ such that for all $s,\tau\in\left[  t,t+\bar{h}\right]
$, $\mathbb{P}$-a.s.
\[
\left\vert U\left(  s\right)  -U\left(  \tau\right)  \right\vert \leq
D_{R}\left\vert s-\tau\right\vert ^{1-\alpha},\quad\left\vert
V\left(  s\right)  -V\left(  \tau\right)  \right\vert \leq\tilde{D}_{R}\left\vert s-\tau\right\vert
^{min\{\beta,1-\alpha\}}%
\]
and
\[
x+\int_{t}^{s}(u+U(r))dr+ \int_{t}^{s} (v+V(r))dB_{r}^{H} \in K(s)  ,
\]
where the constants $D_{R}$, $\tilde{D}_{R}%
$  depend only on $R$, $L_{R}%
$,$M_{0,T}$, $M_{0}$, $L_{0}$,
$T$, $\alpha$, $\beta$, $\Lambda_{\alpha}\left(
B^{H}\right)  $.
\end{definition}
\textbf{Remark.}

$\bullet$ From~\cite{pi}, we can always assume that
$0<\bar{h}\leq1$.

$\bullet$ The definition of $S_{\varphi(K(t))}(t,\varphi(x))$ is the same to $S_{K(t)}(t,x)$, only changes the condition
\[
x+\int_{t}^{s}(u+U(r))dr+ \int_{t}^{s} (v+V(r))dB_{r}^{H} \in K(s)  ,
\]
\qquad to
\[
\varphi(x)+\int_{t}^{s}(u+U(r))dr+ \int_{t}^{s} (v+V(r))dB_{r}^{H} \in \varphi(K(s)).
\]

Now we recall the main result of~\cite{pi} concerning the stochastic
viability.
\begin{theorem}\label{th2.6}
Let $\mathcal{K}=\left\{  K\left(  t\right)  :t\in\left[  0,T\right]  \right\}
$\textit{, }$K\left(  t\right)  =\overline{K\left(  t\right)  }\subset
\mathbb{R}^{d}$. Assume that $(\mathbf{H}_{1})$ and $(\mathbf{H}_{2})$ are satisfied with $\frac{1}{2}<H<1,~ 1-H<\beta, ~\delta>\frac{1-H}{H}$.
Let $1-H<\alpha<\alpha_{0}$. Then the following assertions are equivalent:

\begin{itemize}
\item[(I)] $\mathcal{K}$ is viable for the fractional SDE, i.e. for
all $t\in\left[  0,T\right]  $ and for all $x\in K\left(  t\right)
$ there exists a solution $X^{t,x}\left(  \omega,\cdot\right)  \in
C^{1-\alpha}\left(  \left[  t,T\right]  ;\mathbb{R}^{d}\right)  $ of
the equation
\[
X_{s}^{t,x}=x+\int_{t}^{s}b(r,X_{r}^{t,x})dr+\int_{t}^{s}%
\sigma(r,X_{r}^{t,x})dB_{r}^{H},\;\;s\in\lbrack t,T],\ a.s.\;\omega\in\Omega,
\]
and
\[
X_{s}^{t,x}\in K\left(  s\right)  ,\quad\quad\forall~s\in\left[  t,T\right]
.
\]

\item[(II)] For all $t\in\left[  0,T\right]  $ and all $x\in
K\left(  t\right)  ,$ $\left(  b\left(  t,x\right)  ,\sigma\left(  t,x\right)  \right)  $\textit{ is
}$\left(  1-\alpha\right)  $-fractional $B^{H}$-contingent to $K\left(  t\right)  $ in $x$ .
\item[(III)] For all $t\in\left[  0,T\right]  $ and all $x\in
K(t)  ,$\textit{ }$\left(  b\left(  t,x\right)  ,\sigma\left(  t,x\right)  \right)  $ is
$\left(  1-\alpha\right)  $-fractional $B^{H}$-tangent to $K(t)  $ in $x$ .
\end{itemize}
\end{theorem}
\textbf{Remark.} The assertion (III) is given only for the
deterministic case in~\cite{pi}. In fact we can obtain the
stochastic case from the deterministic one in the same manner as
that (II) is obtained.
\smallskip

Under the same assumptions in Theorem \ref{th2.6}, it follows:

\begin{corollary}\label{cor2.7}
If $K$ is independent of $t$, the following assertions are equivalent:
\begin{itemize}
\item[(j)] $K$ is viable for the fractional SDE (\ref{equ2.3}).
\item[(jj)] For all $t\in\left[  0,T\right]  $ and all $x\in
\partial K  ,$ $\left(  b\left(  t,x\right)  ,\sigma\left(  t,x\right)  \right)  $ is
$\left(  1-\alpha\right)  $-fractional $B^{H}$-contingent to $K$ in $x$ .
\item[(jjj)] For all $t\in\left[  0,T\right]  $ and all $x\in
\partial K  ,$\textit{ }$\left(  b\left(  t,x\right)  ,\sigma\left(  t,x\right)  \right)  $ is
$\left(  1-\alpha\right)  $-fractional $B^{H}$-tangent to $K$ in $x$ .
\end{itemize}
\end{corollary}
\textbf{Proof.}
When $K$ is independent of  $t$, just using Theorem \ref{th2.6}, it's obvious that $(j)\Rightarrow(jj)\Rightarrow(jjj)$.
Now we only need prove  $(jjj)\Rightarrow(j)$, In fact we will prove $(jjj)\Rightarrow(III)$, and then we will get our result.

Let $t\in[0,T] $ and  $\forall x\in K\setminus \partial K$, Since
$X^{t,x}$ is continuous, then there exists a random variable
$\bar{h}$, such that for all $s\in [t,t+\bar{h}]$,
\begin{equation*}
X_{s}^{t,x}=x+{\displaystyle\int_{t}^{s}}b(r,X_{r}^{t,x})dr+{\displaystyle\int_{t}^{s}}\sigma(r,X_{r}^{t,x})dB_{r}^{H}\in K.
\end{equation*}
we have for all $s\in [t,t+\bar{h}]$,
\begin{equation*}
X_{s}^{t,x}=x+{\displaystyle\int_{t}^{s}}[b(t,x)+U(r)]dr+{\displaystyle\int_{t}^{s}}[\sigma(t,x)+V(r)]dB_{r}^{H}\in K.
\end{equation*}
where
\[U(r)=b(r,X_{r}^{t,x})-b(t,x), \quad  V(r)=\sigma(r,X_{r}^{t,x})-\sigma(t,x)
\]
clearly that $\left(  b\left(  t,x\right)  ,\sigma\left(  t,x\right)
\right)  $ is $\left(  1-\alpha\right)  $-fractional $B^{H}$-tangent
to $K  $ in $x$. Together with $(jjj)$, we have that for all
$t\in\left[ 0,T\right]  $ and all $x\in K  ,$\textit{ }$\left(
b\left( t,x\right)  ,\sigma\left(  t,x\right)  \right)  $ is $\left(
1-\alpha\right)  $-fractional $B^{H}$-tangent to $K$ in $x$. This is
just $(III)$ for the case that $K$ is independent of $t$.
\begin{flushright}
$\Box$
\end{flushright}
\section{Results and Applications}\label{sec3}

The next two theorems are our main theorems, firstly we extend \textit{Stochastic Tangent Sets to Direct Images}
which is introduced by J.P.Aubin, and G.Da Prato~\cite{pc} (1990) to the fBM framework.

\begin{theorem}\label{th3.1}
Assume that $(\mathbf{H}_{1})$ and $(\mathbf{H}_{2})$ are satisfied. Let
%$\mathcal{K}=\left\{  K\left(  t\right)  :t\in\left[  0,T\right]  \right\}
$K\left(  t\right)  =\overline{K\left(  t\right)  }\subset
\mathbb{R}^{d}, t\in\left[  0,T\right]$ and $S_{K(t)}(t,x)$ the $\left(  1-\alpha\right)  $-fractional $B^{H}$-tangent set to $K$ in $x$.
Let $\varphi$ be a $C^{2}$ map from $\mathbb{R}^d$ to $\mathbb{R}^m$ with a bounded second derivative.
If
\[(b(t,x),\sigma(t,x))\in S_{K(t)}(t,x)
\]
then
\[(\varphi^{\prime}(x)b(t,x),\varphi^{\prime}(x)\sigma(t,x))\in S_{\varphi(K(t))}(t,\varphi(x)).
\]
\end{theorem}
\medskip

Also we can prove the \textit{ Stochastic Tangent Sets to Inverse Images} in the fBM form.\\

We introduce a space $\mathcal{H}$ of the functions $\varphi:
\mathbb{R}^{d}\rightarrow\mathbb{R}^{m}$ of class $C^{2}$, with a
bounded and Lipschitz continuous second derivative and there exist
$a_{\varphi}<b_{\varphi}$ and some constants $M>0$, $L>0$ such that
for all $a_{\varphi}\leq|x|\leq b_{\varphi}$, the matrix
$\varphi^{\prime}(x)$ has a right inverse denoted by
$\varphi^{\prime}(x)^{+}$ satisfying
\[
\begin{array}
[c]{ll}
(1)  \; \qquad \left|[\varphi^{\prime}(x)^{+}]^{\prime}\right|\leq M,\\
(2)  \;
\qquad\left|[\varphi^{\prime}(x)^{+}]^{\prime}-[\varphi^{\prime}(y)^{+}]^{\prime}\right|\leq
L|x-y|.
\end{array}
\]
\begin{theorem}\label{th3.2}
Assume that $(\mathbf{H}_{1})$ and $(\mathbf{H}_{2})$ are satisfied.
Let $K\left(  t\right)  =\overline{K\left(  t\right)  }\subset
\mathbb{R}^{d}, t\in\left[  0,T\right]$ and
$\varphi\in \mathcal{H}$,
\textit{then for} every $\varepsilon>0$ and $a_{\varphi}+\varepsilon\leq|x|\leq b_{\varphi}-\varepsilon$, then
\[(b(t,x),\sigma(t,x))\in S_{\varphi^{-1}(\varphi(K(t)))}(t,x)
\]
if and only if
\[(\varphi^{\prime}(x)b(t,x),\varphi^{\prime}(x)\sigma(t,x))\in S_{\varphi(K(t))}(t,\varphi(x)).
\]
\end{theorem}
\medskip

Using Theorem \ref{th3.2}, we can get the deterministic sufficient
and necessary conditions for viability when $K$ takes some particular forms. Firstly we give some Lemmas.
\begin{lemma}\label{lem3.3}
Let $K$ be the unit sphere, then for all $x\in K$, $(b(t,x),\sigma(t,x))\in S_{K}(t,x)$ if and only if
\[\langle x,b(t,x)\rangle=0, \qquad \langle x,\sigma(t,x)\rangle=0.
\]
\end{lemma}
\begin{lemma}\label{lem3.4}
Let $K=\{x\in \mathbb{R}^{d}; r\leq|x|\leq R\}$ then for all $x$, such that $|x|=R$, $(b(t,x),\sigma(t,x))\in S_{K}(t,x)$ if and only if
\[\langle x,b(t,x)\rangle\leq 0, \qquad \langle x,\sigma(t,x)\rangle=0
\]
and for all $x$, such that $|x|=r$, $(b(t,x),\sigma(t,x))\in S_{K}(t,x)$ if and only if
\[\langle x,b(t,x)\rangle\ge 0, \qquad \langle x,\sigma(t,x)\rangle=0.
\]
\end{lemma}
\begin{lemma}\label{lem3.5}
Let $K$ be the unit ball, then for all $x$, such that $|x|=1$, $(b(t,x),\sigma(t,x))\in S_{K}(t,x)$ if and only if
\[\langle x,b(t,x)\rangle\leq 0, \qquad \langle x,\sigma(t,x)\rangle=0.
\]
\end{lemma}
\medskip

Just as Corollary \ref{cor2.7} said, considering that if we want to get the conditions for the viability of $K$,
we only need to think about the starting point $x\in\partial K$. Then together with Lemma \ref{lem3.3} and \ref{lem3.5}, it is obviously that

\begin{proposition}
Let $(\mathbf{H}_{1})$, $(\mathbf{H}_{2})$ be satisfied, $1-H<\alpha<\alpha_{0}$ and $K$ is the unit sphere.
Then the following assertions are equivalent:
\begin{itemize}
\item[(I)] $K$ is viable for the fractional SDE (\ref{equ2.3}).
\item[(II)] For all $t\in\left[  0,T\right]  $ and all $x\in
K$  ,
\[
\langle x,b(t,x)\rangle=0,\quad \langle x,\sigma(t,x)\rangle=0.
\]
\end{itemize}
\end{proposition}
\begin{proposition}\label{pro}
Let $(\mathbf{H}_{1})$, $(\mathbf{H}_{2})$ be satisfied, $1-H<\alpha<\alpha_{0}$ and $K$ is the unit ball.
Then the following assertions are equivalent:
\begin{itemize}
\item[(I)] $K$ is viable for the fractional SDE (\ref{equ2.3}).
\item[(II)] For all $t\in\left[  0,T\right]  $ and all $|x|=1$,
\[
\langle x,b(t,x)\rangle\leq 0,\quad \langle x,\sigma(t,x)\rangle=0.
\]
\end{itemize}
\end{proposition}
\begin{corollary}\label{cor3.8}
Consider the SDE on $\mathbb{R}$,
\[
X_{s}=x+\int_{t}^{s}b(r,X_{r})dr+\int_{t}^{s}\sigma\left(  r,X_{r}\right)
dB_{r}^{H},\,\;s\in\left[  t,T\right]  .
\]
$B^{H}=\left\{  B^{H}_{t},t\geq0\right\}  $ is a
fractional Brownian motion. $b, \sigma$ satisfy the assumptions $(\mathbf{H}_{1}), (\mathbf{H}_{2})$. Then for any $t\in[0,T]$ and
every $x\ge0$ the equation has a positive solution if and only if
\[b(t,0)\geq0, \quad\sigma(t,0)=0, \quad\forall t\in[0,T].
\]

\end{corollary}

\textbf{Proof.}
In fact we take $K=[0,+\infty)$, the problem is just that $K$ is viable for the fractional SDE.
We can use $x=\tan\frac{\pi}{4}(y+1)$ and we get $y=\frac{4}{\pi}\arctan x-1$,
it just maps $[0,+\infty)$ to $[-1,1]$, and  using Proposition \ref{pro} and It\^{o} formula of fractional SDE (see~\cite{pk}), we have
\[b(t,0)\geq0, \quad\sigma(t,0)=0, \quad\forall t\in[0,T].
\]
\begin{flushright}
$\Box$
\end{flushright}
The most interesting application is the characterization of comparison theorem.\\
Let us firstly consider the linear case.
\begin{corollary}
Consider the linear two dimensional decoupled system
\[
\left\{
\begin{array}
[c]{l}%
X_{s}^{t,x}=x+{\displaystyle\int_{t}^{s}}(f(r)X_{r}^{t,x}+f_{1}(r))dr+{\displaystyle\int_{t}^{s}}(g(r)X_{r}^{t,x}+g_{1}(r))dB^{H}_{r},\smallskip\;
s\in[t,T] \medskip\\
Y_{s}^{t,y}=y+{\displaystyle\int_{t}^{s}}(f(r)Y_{r}^{t,y}+f_{2}(r))dr+{\displaystyle\int_{t}^{s}}%
(g(r)Y_{r}^{t,y}+g_{2}(r))dB^{H}_{r},\smallskip\;s\in[t,T]
\end{array}
\right.
\]
then
\[\textit{for any}~t\in[0,T] ~\textit{and every}~x\leq y, ~ X_{s}^{t,x}\leq Y_{s}^{t,y},~\forall s\in[t,T].\]
\[\Longleftrightarrow ~f_{1}(t)\leq f_{2}(t),~ g_{1}(t)=g_{2}(t),~\forall t\in[0,T] .
\]
\end{corollary}

\textbf{Proof.} In fact we set
$Z_{s}^{t,z}=Y_{s}^{t,y}-X_{s}^{t,x}$, where $z=y-x\ge0$, then we
can change the problem to $Z_{s}^{t,z}\ge0$, it means that  for any
$t\in[0,T]$ and every $z\ge0$ the fractional SDE of $Z_{s}^{t,z}$
has a positive solution. Then using Corollary \ref{cor3.8}, we can
easily prove this corollary.
\begin{flushright}
$\Box$
\end{flushright}
In general, we have
\begin{theorem}\label{th3.10}(Comparison theorem)
Consider the two dimensional decoupled system%
\[
\left\{
\begin{array}
[c]{l}%
X_{s}^{t,x}=x+{\displaystyle\int_{t}^{s}}(b_{1}(r,X_{r}^{t,x}))dr+{\displaystyle\int_{t}^{s}}(\sigma_{1}(r,X_{r}^{t,x}))dB^{H}_{r},
\smallskip\;s\in[t,T] \medskip\\
Y_{s}^{t,y}=y+{\displaystyle\int_{t}^{s}}(b_{2}(r,Y_{r}^{t,y}))dr+{\displaystyle\int_{t}^{s}}(\sigma_{2}(r,Y_{r}^{t,y}))dB^{H}_{r},
\smallskip\;s\in[t,T]
\end{array}
\right.
\]
then
\[\textit{for any}~ t\in[0,T] ~\textit{and every}~x\leq y,~X_{s}^{t,x}\leq Y_{s}^{t,y},~\forall s\in[t,T]
\]
\[\Longleftrightarrow b_{1}(t,z)\leq b_{2}(t,z),~ \sigma_{1}(t,z)=\sigma_{2}(t,z),~\forall t\in[0,T], ~\forall z\in\mathbb{R}.
\]
\end{theorem}

we will give the proof of this result in Section \ref{sec5}.

\section{Proofs of main results}
This section is devoted to the proofs of the main results which have been given in Section \ref{sec3}.

Firstly we present some auxiliary Lemmas which will be used in the
sequel.
\subsection{Auxiliary Results}

\begin{lemma}\label{lem4.1}
Given two stochastic process
\[
\begin{array}
[c]{ll}
U=U^{t,x}:& \Omega\times\left[  t,t+\bar{h}\right]  \rightarrow\mathbb{R}^{d} ,~\ U(t)=0
\\
V=V^{t,x}:& \Omega\times\left[  t,t+\bar{h}\right]  \rightarrow\mathbb{R}^{d} ,~\ V(t)=0
\end{array}
\]
such that for all $s,\tau\in\left[  t,t+\bar{h}\right]  $ and
for every $R>0$ with $\left\vert x\right\vert \leq R:$
\[
\begin{array}
[c]{ll}
\left\vert U\left(  s\right)  -U\left(  \tau\right)  \right\vert & \leq
D_{R}\left\vert s-\tau\right\vert ^{1-\alpha},\\
\left\vert
V\left(  s\right)  -V\left(  \tau\right)  \right\vert & \leq\tilde{D}_{R}\left\vert s-\tau\right\vert
^{min\{\beta,1-\alpha\}}.
\end{array}
\]
then for all $t\leq\tau\leq s\leq t+\bar{h}$,
\[
\begin{array}
[c]{ll}
\left(  a\right)  \; & \left\vert
%TCIMACRO{\dint _{t}^{s}}%
%BeginExpansion
{\displaystyle\int_{\tau}^{s}}
%EndExpansion
U(r)dr\right\vert \leq
D_{R}\left(  s-t\right)  ^{1-\alpha}(s-\tau)\\
\left(  b\right)  \; & \left\vert
%TCIMACRO{\dint _{t}^{s}}%
%BeginExpansion
{\displaystyle\int_{\tau}^{s}}
%EndExpansion
V(r)  dB^{H}_{r}\right\vert \leq
C_{R}(\alpha,\beta)\tilde{D}_{R}\Lambda_{\alpha}(B^{H})(s-t)^{\min\{\beta,1-\alpha\}}\left(
s-\tau\right)  ^{1-\alpha}.
\end{array}
\]
where $C_{R}(\alpha,\beta)$ depends only on $R$, $\alpha$, and
$\beta$.
\end{lemma}
\textbf{Proof.}

$(a)$ we have
\begin{eqnarray*}
\left\vert
%TCIMACRO{\dint _{t}^{s}}%
%BeginExpansion
{\displaystyle\int_{\tau}^{s}}
%EndExpansion
U(r) dr\right\vert
& = &
\left\vert
%TCIMACRO{\dint _{t}^{s}}%
%BeginExpansion
{\displaystyle\int_{\tau}^{s}}
%EndExpansion
[U(r)-U(t)]  dr\right\vert
{}
\nonumber\\
& \leq & {} D_{R}\left\vert
%TCIMACRO{\dint _{t}^{s}}%
%BeginExpansion
{\displaystyle\int_{\tau}^{s}}
%EndExpansion
(r-t)^{1-\alpha}  dr\right\vert {}
\nonumber\\
& \leq & {} D_{R}(s-t)^{1-\alpha}(s-\tau).
\end{eqnarray*}

$(d)$
\begin{eqnarray*}
\left\vert
%TCIMACRO{\dint _{t}^{s}}%
%BeginExpansion
{\displaystyle\int_{\tau}^{s}}
%EndExpansion
V(r) dB_{r}^{H}\right\vert
& = &
\left\vert
%TCIMACRO{\dint _{t}^{s}}%
%BeginExpansion
{\displaystyle\int_{\tau}^{s}}
%EndExpansion
[V(r)-V(t)]  dB_{r}^{H}\right\vert
{}
\nonumber\\
& \leq & {}
\Lambda_{\alpha}(B^{H})\|V\|_{\alpha,1;[\tau,s]}
{}
\nonumber\\
& \leq & {}
\Lambda_{\alpha}(B^{H})\int_{\tau}^{s}\left[  \frac{|V(r)-V(t)|}{\left(  r-\tau\right)  ^{\alpha}}+\int_{\tau}%
^{r}\frac{|V(r)-V(y)|}{(r-y)^{\alpha+1}}dy\right] dr {}
\nonumber\\
& \leq & {}
\tilde{D}_{R}\Lambda_{\alpha}(B^{H})\int_{\tau}^{s}\left[  \frac{(r-t)^{\min\{\beta,1-\alpha\}}}{\left(  r-\tau\right)  ^{\alpha}}+\int_{\tau}%
^{r}\frac{(r-y)^{\min\{\beta,1-\alpha\}}}{(r-y)^{\alpha+1}}dy\right]
dr {}
\nonumber\\
& \leq & {}
\tilde{D}_{R}\Lambda_{\alpha}(B^{H})\Big[\frac{1}{1-\alpha}(s-t)^{\min\{\beta,1-\alpha\}}\left(
s-\tau\right)  ^{1-\alpha}
 {}
\nonumber\\
&   & {}
\qquad\qquad\qquad\qquad+\int_{\tau}^{s}\int_{\tau}%
^{r}(r-y)^{\min\{\beta-\alpha,1-2\alpha\}-1}dy dr {}\Big]
\nonumber\\
& \leq & {}
C_{R}(\alpha,\beta)\tilde{D}_{R}\Lambda_{\alpha}(B^{H})(s-t)^{\min\{\beta,1-\alpha\}}\left(
s-\tau\right)  ^{1-\alpha}.
\end{eqnarray*}
\begin{flushright}
$\Box$
\end{flushright}
\textbf{Remark.} From $(a)$ and $(b)$, just taking $\tau=t$, we have
\[
\begin{array}
[c]{ll} \left(  a^{\prime}\right)  \; & \left\vert
%TCIMACRO{\dint _{t}^{s}}%
%BeginExpansion
{\displaystyle\int_{t}^{s}}
%EndExpansion
U(r)dr\right\vert \leq
D_{R}\left(  s-t\right)  ^{2-\alpha}\\
\left(  b^{\prime}\right)  \; & \left\vert
%TCIMACRO{\dint _{t}^{s}}%
%BeginExpansion
{\displaystyle\int_{t}^{s}}
%EndExpansion
V(r)  dB^{H}_{r}\right\vert \leq
C_{R}(\alpha,\beta)\tilde{D}_{R}\Lambda_{\alpha}(B^{H})(s-t)^{1+\min\{\beta-\alpha,1-2\alpha\}}.
\end{array}
\]
\begin{lemma}\label{lem4.2}
Given two stochastic process $U=U^{t,x}, V=V^{t,x}$ which satisfy the conditions in Lemma \ref{lem4.1}, and $\varphi\in\mathcal{H}$, let
\begin{eqnarray*}
f(r,y) & = &  \varphi^{\prime}(y)^{+}\left[U(r)-(\varphi^{\prime}(y)-\varphi^{\prime}(x))b(t,x)\right]\\
g(r,y) & = &  \varphi^{\prime}(y)^{+}\left[V(r)-(\varphi^{\prime}(y)-\varphi^{\prime}(x))\sigma(t,x)\right]
\end{eqnarray*}
then for $\alpha\in\left(  1-H,\alpha_{0}\right)$ and for every
$\delta_{0}>0$ there exists a random variable
$\bar{h}_{1}=\bar{h}_{1}^{t,x}$ such that for
$a_{\varphi}+2\delta_{0}\leq|x|\leq b_{\varphi}-2\delta_{0}$ and
$\mathbb{P}$-a.s. $\omega\in\Omega$, the following SDE
\[
\xi_{s}=x+\int_{t}^{s}(b(t,x)+f(r,\xi_{r}))dr+\int_{t}^{s}(\sigma(t,x)+g(r,\xi_{r}))dB^{H}(r),
~s\in[t,t+\bar{h}_{1}],
\]
has a unique solution $\xi_{\cdot}\left(  \omega\right) \in
L^{0}\left(  \Omega,\mathcal{F},\mathbb{P\,};
W^{\alpha,\infty}(t,T;\mathbb{R}^{d})\right).$ \\Moreover
$\mathbb{P}$-a.s. $\xi_{\cdot}\left(  \omega\right) \in
C^{1-\alpha}\left(  t,t+\bar{h}_{1};\mathbb{R}^{d}\right)$.
\end{lemma}
\textbf{Proof.} From \cite{po} Theorem(the partition of unity) p.61,
we have that for every $\delta_{0}>0$, there exists one function
$\alpha(x)\in C^{\infty}(\mathbb{R}^{d})$ such that $\alpha(x)=1$
for $a_{\varphi}+\delta_{0}\leq|x|\leq b_{\varphi}-\delta_{0}$ and
$\alpha(x)=0$ for $|x|\ge b_{\varphi}~ or ~|x|\leq a_{\varphi}$,
then we define
\[\tilde{f}(t,y)=\alpha(y)f(t,y)=\left\{
\begin{array}
[c]{ll}%
f(t,y), \quad  & a_{\varphi}+\delta_{0}\leq|y|\leq b_{\varphi}-\delta_{0}\\
\alpha(y)f(t,y) \quad & a_{\varphi}\leq|y|\leq a_{\varphi}+\delta_{0}, ~or ~b_{\varphi}-\delta_{0}\leq|y|\leq b_{\varphi}\\
0, \quad & |y|\ge b_{\varphi}, ~or~|y|\leq a_{\varphi}.
\end{array}
\right.
\]
and we define $\tilde{g}(t,y)$ in the same method and then we consider the following SDE
\begin{equation}\label{equ1}
\tilde{\xi}_{s}=x+\int_{t}^{s}(b(t,x)+\tilde{f}(r,\tilde{\xi}_{r}))dr+\int_{t}^{s}(\sigma(t,x)+\tilde{g}(r,\tilde{\xi}_{r}))dB^{H}(r),
~s\in[t,t+\bar{h}].
\end{equation}
Since $\varphi\in\mathcal{H}$ and $U=U^{t,x}, V=V^{t,x}$ satisfy the
conditions in Lemma \ref{lem4.1}, we can verify that for
$\alpha\in\left(  1-H,\alpha_{0}\right),$ $\tilde{f}(t,y)$, and
$\tilde{g}(t,y)$ satisfy the conditions in
($\mathbf{H}_{1}$),($\mathbf{H}_{2}$) in ~\cite{pi} where the
constants $M_{0}, M_{R},L_{0}, L_{R}$ depend on $\omega$, then the
SDE (\ref{equ1}) has a unique solution $\tilde{\xi}_{\cdot}\left(
\omega\right) \in L^{0}\left(
\Omega,\mathcal{F},\mathbb{P\,};W^{\alpha,\infty}(t,T;\mathbb{R}^{d})\right)$
for all $\alpha\in\left(  1-H,\alpha_{0}\right)  .$ And moreover
$\mathbb{P}$-a.s. $\tilde{\xi}_{\cdot}\left(  \omega\right) \in
C^{1-\alpha}\left(  t,t+\bar{h};\mathbb{R}^{d}\right)$. Since
$a_{\varphi}+2\delta_{0}\leq|x|\leq b_{\varphi}-2\delta_{0}$ then
there exists a random variable $\bar{h}_{1}=\bar{h}_{1}^{t,x}$, such
that $\mathbb{P}$-a.s. $a_{\varphi}+\delta_{0}\leq|\tilde{\xi}|\leq
b_{\varphi}-\delta_{0}$, then for $s\in[t,t+\bar{h}_{1}]$, the SDE
(\ref{equ1}) becomes $\mathbb{P}$-a.s.
\begin{equation*}
\tilde{\xi}_{s}=x+\int_{t}^{s}(b(t,x)+f(r,\tilde{\xi}_{r}))dr+\int_{t}^{s}(\sigma(t,x)+g(r,\tilde{\xi}_{r}))dB^{H}(r),
~s\in[t,t+\bar{h}_{1}].
\end{equation*}
just taking $\xi_{s}=\tilde{\xi}_{s}, s\in[t,t+\bar{h}_{1}]$, and together with the uniqueness of $\tilde{\xi}_{s}$,  then we finish our proof.
\begin{flushright}
$\Box$
\end{flushright}

\subsection{Proof of Theorem 3.1 and Theorem 3.2}

\textbf{Proof of Theorem \ref{th3.1}}

Since $(b(t,x),\sigma(t,x))\in S_{K(t)}(t,x)$, then
there exist a random variable $\bar{h}=\bar{h}^{t,x}>0,$ and two stochastic process
\[
\begin{array}
[c]{ll}
U=U^{t,x}: & \Omega\times\left[  t,t+\bar{h}\right]  \rightarrow\mathbb{R}^{d} ,~\ U(t)=0
\\
V=V^{t,x}:& \Omega\times\left[  t,t+\bar{h}\right]  \rightarrow\mathbb{R}^{d} ,~\ V(t)=0
\end{array}
\]
such that for all $s,\tau\in\left[  t,t+\bar{h}\right]  $ and for
every $R>0$ and $\left\vert x\right\vert \leq R:$
\[
\left\vert U\left(  s\right)  -U\left(  \tau\right)  \right\vert \leq
D_{R}\left\vert s-\tau\right\vert ^{1-\alpha},\quad\left\vert
V\left(  s\right)  -V\left(  \tau\right)  \right\vert \leq\tilde{D}_{R}\left\vert s-\tau\right\vert
^{min\{\beta,1-\alpha\}}%
\]
and
\[
x+\int_{t}^{s}(b(t,x)+U(r))dr+ \int_{t}^{s} (\sigma(t,x)+V(r))dB^{H}(r) \in K(s)  ,
\]
where  $D_{R}$, $\tilde{D}_{R}$, depend only on $R$, $L_{R}%
$, $M_{0,T}$, $M_{0}$, $L_{0}$,
$T$, $\alpha$, $\beta$, $\Lambda_{\alpha}(B^{H})$.

Let
\[\eta_{s}=x+\int_{t}^{s}(b(t,x)+U(r))dr+\int_{t}^{s}(\sigma(t,x)+V(r))dB^{H}(r)
\]
and from Lemma \ref{lem4.1} and $H-\epsilon$ H\"{o}lder continuous property of fractional Brownian motion,
it follows that for all $s,\tau\in\left[  t,t+\bar{h}\right]  $,
\[|\eta_{s}-\eta_{\tau}|\leq \zeta(s-\tau)^{1-\alpha}.
\]
According to the fractional It\^{o} formula (see Yuliya
S.Mishura~\cite{pk}), We have for all $s\in\left[ t,t+\bar{h}\right]
$
\[
\varphi\Big(x+\int_{t}^{s}(b(t,x)+U(r))dr+\int_{t}^{s}(\sigma(t,x)+V(r))dB^{H}(r)\Big)
\]
\[
=\varphi(x)+\int_{t}^{s}[\varphi^{\prime}(\eta_{r})(b(t,x)+U(r))]dr+\int_{t}^{s}[\varphi^{\prime}(\eta_{r})(\sigma(t,x)+V(r))]dB^{H}(r)
\]
\[
=\varphi(x)+\int_{t}^{s}[\varphi^{\prime}(x)b(t,x)+U_{1}(r)]dr+\int_{t}^{s}[\varphi^{\prime}(x)\sigma(t,x)+V_{1}(r)]dB^{H}(r)
\]
where
\begin{eqnarray*}
U_{1}(r) & = & \varphi^{\prime}(\eta_{r})U(r)+(\varphi^{\prime}(\eta_{r})-\varphi^{\prime}(x))b(t,x)\\
V_{1}(r) & = & \varphi^{\prime}(\eta_{r})V(r)+(\varphi^{\prime}(\eta_{r})-\varphi^{\prime}(x))\sigma(t,x)
\end{eqnarray*}
Then
\[
\varphi(x)+\int_{t}^{s}[\varphi^{\prime}(x)b(t,x)+U_{1}(r)]dr+\int_{t}^{s}[\varphi^{\prime}(x)\sigma(t,x)+V_{1}(r)]dB^{H}(r)
\]
\[=
\varphi\Big(x+\int_{t}^{s}(b(t,x)+U(r))dr+\int_{t}^{s}(\sigma(t,x)+V(r))dB^{H}(r)\Big)\in\varphi(K(s))
\]
and it's easy to verify that
\[U_{1}(t)=0, \quad V_{1}(t)=0
\]
For all $s,\tau\in\left[  t,t+\bar{h}\right]  $ and for every $R>0$
and $\left\vert x\right\vert \leq R,$ Using the Lipschitz continuity
of $\varphi^{\prime}$ and ($\mathbf{H}_{2}$), we obtain that
\begin{eqnarray*}
|U_{1}(s)-U_{1}(\tau)|&\leq& |\varphi^{\prime}(\eta_{\tau})||U(s)-U(\tau)|+(|U(s)|+|b(t,x)|)|\varphi^{\prime}(\eta_{s})-\varphi^{\prime}(\eta_{\tau})|\\
&\leq& \theta_{1}|s-\tau|^{1-\alpha}+ \theta_{2}|\eta_{s}-\eta_{\tau}|\\
&\leq& \theta|s-\tau|^{1-\alpha}
\end{eqnarray*}
Similarly we can proof that
\begin{eqnarray*}
|V_{1}(s)-V_{1}(\tau)|\leq
\tilde{\theta}|s-\tau|^{min\{\beta,1-\alpha\}}
\end{eqnarray*}
The H\"{o}lder constants $\theta$, $\tilde{\theta}$ are random variables which
depend only on $R$, $L_{R}$,  $M_{0}$, $M_{0}$, $L_{0}$,
$T$, $\alpha$, $\beta$, $\Lambda_{\alpha}\left(
B^{H}\right)  $.

This means that
\[(\varphi^{\prime}(x)b(t,x),\varphi^{\prime}(x)\sigma(t,x))\in S_{\varphi(K(t))}(t,\varphi(x)).
\]
\begin{flushright}
$\Box$
\end{flushright}
\textbf{Proof of Theorem \ref{th3.2}}

We shall only have to prove that from
$(\varphi^{\prime}(x)b(t,x),\varphi^{\prime}(x)\sigma(t,x))\in S_{\varphi(K(t))}(t,\varphi(x)),$
we infer that $(b(t,x),\sigma(t,x))\in S_{\varphi^{-1}(\varphi(K(t)))}(t,x)$.

Since
\[(\varphi^{\prime}(x)b(t,x),\varphi^{\prime}(x)\sigma(t,x))\in S_{\varphi(K(t))}(t,\varphi(x)).
\]
then for $x\in K(t)$, there exist  a random variable
$\bar{h}=\bar{h}^{t,x}>0$ and two stochastic process,
\[
\begin{array}
[c]{ll} U_{1}=U_{1}^{t,x}:& \Omega\times\left[  t,t+\bar{h}\right]
\rightarrow\mathbb{R}^{m} ,~\ U_{1}(t)=0,
\\
V_{1}=V_{1}^{t,x}:& \Omega\times\left[  t,t+\bar{h}\right]
\rightarrow\mathbb{R}^{m} ,~\ V_{1}(t)=0
\end{array}
\]
such that for all $s,\tau\in\left[  t,t+\bar{h}\right]  $ and for every $R>0$ and $\left\vert x\right\vert \leq R,$

\[
\left\vert U_{1}\left(  s\right)  -U_{1}\left(  \tau\right)  \right\vert \leq
D_{R}\left\vert s-\tau\right\vert ^{1-\alpha},\quad\left\vert
V_{1}\left(  s\right)  -V_{1}\left(  \tau\right)  \right\vert \leq\tilde{D}_{R}\left\vert s-\tau\right\vert
^{min\{\beta,1-\alpha\}}%
\]
and
\[
\varphi(x)+\int_{t}^{s}(\varphi^{\prime}(x)b(t,x)+U_{1}(r))dr+ \int_{t}^{s} (\varphi^{\prime}(x)\sigma(t,x)+V_{1}(r))dB^{H}(r) \in \varphi(K(s)).
\]
Let
\begin{eqnarray*}
f(r,y) & = &  \varphi^{\prime}(y)^{+}\left[U_{1}(r)-(\varphi^{\prime}(y)-\varphi^{\prime}(x))b(t,x)\right],\\
g(r,y) & = &
\varphi^{\prime}(y)^{+}\left[V_{1}(r)-(\varphi^{\prime}(y)-\varphi^{\prime}(x))\sigma(t,x)\right],
\end{eqnarray*}
where $\varphi^{\prime}(y)^{+}$ is the right inverse of
$\varphi^{\prime}(y)$. By Lemma \ref{lem4.2}, for every
$\delta_{0}>0$ and $a_{\varphi}+2\delta_{0}\leq|x|\leq
b_{\varphi}-2\delta_{0}$,  there exists a random variable
$\bar{h}_{1}$ such that for $\mathbb{P}$-a.s. $\omega\in\Omega$ the
following SDE
\[
\xi_{s}=x+\int_{t}^{s}(b(t,x)+f(r,\xi_{r}))dr+\int_{t}^{s}(\sigma(t,x)+g(r,\xi_{r}))dB^{H}(r),
~s\in[t,t+\bar{h}_{1}],
\]
has a unique solution $\xi_{\cdot}\left(  \omega\right)$. Then with
\begin{eqnarray*}
U(r) & = &  \varphi^{\prime}(\xi_{r})^{+}\left[U_{1}(r)-(\varphi^{\prime}(\xi_{r})-\varphi^{\prime}(x))b(t,x)\right] \quad\text{and}\\
V(r) & = & \varphi^{\prime}(\xi_{r})^{+}\left[V_{1}(r)-(\varphi^{\prime}(\xi_{r})-\varphi^{\prime}(x))\sigma(t,x)\right]
\end{eqnarray*}
according to the fractional It\^{o} formula, we have for all
$s\in\left[ t,t+\bar{h}_{1}\right]  $
\[
\varphi\Big(x+\int_{t}^{s}(b(t,x)+U(r))dr+\int_{t}^{s}(\sigma(t,x)+V(r))dB^{H}(r)\Big)
\]
\[=
\varphi(x)+\int_{t}^{s}(\varphi^{\prime}(x)b(t,x)+U_{1}(r))dr+\int_{t}^{s}(\varphi^{\prime}(x)\sigma(t,x)+V_{1}(r))dB^{H}(r)\in\varphi(K(s)).
\]
Clearly that
\[U(t)=0,\quad V(t)=0.
\]
Since $\varphi\in\mathcal{H}$ and $\xi_{\cdot}\left(  \omega\right)
\in C^{1-\alpha}\left(  t,t+\bar{h}_{1};\mathbb{R}^{d}\right)$ and
together with ($\mathbf{H}_{1}$) and ($\mathbf{H}_{2}$), it easily
follows
\[
\left\vert U\left(  s\right)  -U\left(  \tau\right)  \right\vert \leq
\theta\left\vert s-\tau\right\vert ^{1-\alpha},\quad\left\vert
V\left(  s\right)  -V\left(  \tau\right)  \right\vert \leq\tilde{\theta}\left\vert s-\tau\right\vert
^{min\{\beta,1-\alpha\}}.%
\]
Then it means that
\[(b(t,x),\sigma(t,x))\in S_{\varphi^{-1}(\varphi(K(t)))}(t,x).
\]
\begin{flushright}
$\Box$
\end{flushright}
\subsection{Proof of Lemmas 3.3, 3.4, 3.5}.
\textbf{Proof of Lemma \ref{lem3.3}}

Firstly, we take $\varphi(x)=|x|^2$, and it's easy to verify that for $\frac{1}{4}\leq|x|\leq4$,
\[\varphi^{\prime}(x)^{+}=\frac{x}{2|x|^2}.
\]
and we can verify that $\varphi\in\mathcal{H}$ taking
$a_{\varphi}=\frac{1}{4}$, $b_{\varphi}=4$,
$\varepsilon=\frac{1}{4}$, then for $x\in K$ we have
$\frac{1}{2}\leq|x|=1\leq4-\frac{1}{4}$, by Theorem \ref{th3.2} we
have
\[(b(t,x),\sigma(t,x))\in S_{K}(t,x)\Leftrightarrow (\langle 2x,b(t,x)\rangle,2x^{\ast}\sigma(t,x))\in S_{1}(t,x^2)
\]

So now it's equivalent to prove
\[(\langle 2x,b(t,x)\rangle,\langle 2x,\sigma(t,x)\rangle)\in S_{1}(t,|x|^2)\Leftrightarrow \langle x,b(t,x)\rangle=0,
\quad \langle x,\sigma(t,x)\rangle=0
\]

\textit{Sufficient.} If $\langle x,b(t,x)\rangle=0,~\langle
x,\sigma(t,x)\rangle=0$, we can take $U(r)\equiv 0, ~V(r)\equiv 0$,
and we have $\forall s\in[t,t+\bar{h}]$ and $|x|=1$
\[|x|^2+\int_{t}^{s}(\langle 2x,b(t,x)\rangle+U(r))dr+\int_{t}^{s}(\langle 2x,\sigma(t,x)\rangle)+V(r))dB^{H}(r)=1,
\]
This means that $(\langle 2x,b(t,x)\rangle, \langle 2x,\sigma(t,x)\rangle)\in S_{1}(t,|x|^2)$.

\textit{Necessary.}  Since $(\langle 2x,b(t,x)\rangle, \langle
2x,\sigma(t,x)\rangle)\in S_{1}(t,|x|^2)$. then there exist a random
variable $\bar{h}=\bar{h}^{t,x}>0,$ and two stochastic process
\[
\begin{array}[c]{ll}
U=U^{t,x}:& \Omega\times\left[  t,t+\bar{h}\right]  \rightarrow\mathbb{R} ,~\ U(t)=0
\\
V=V^{t,x}:& \Omega\times\left[  t,t+\bar{h}\right]  \rightarrow\mathbb{R} ,~\ V(t)=0
\end{array}
\]
such that for all $s,\tau\in\left[  t,t+\bar{h}\right]  $ and for
every $R>0$, $\left\vert x\right\vert \leq R:$
\[
\left\vert U\left(  s\right)  -U\left(  \tau\right)  \right\vert \leq
D_{R}\left\vert s-\tau\right\vert ^{1-\alpha},\quad\left\vert
V\left(  s\right)  -V\left(  \tau\right)  \right\vert \leq\tilde{D}_{R}\left\vert s-\tau\right\vert
^{min\{\beta,1-\alpha\}}%
\]
and
\begin{equation}\label{equ4.1}
|x|^2+\int_{t}^{s}(\langle 2x,b(t,x)\rangle+U(r))dr+\int_{t}^{s}(\langle 2x,\sigma(t,x)\rangle+V(r))dB^{H}(r)=1, ~\mathbb{P}-a.s.
\end{equation}
where  $D_{R}$, $\tilde{D}_{R}%
$,  depend only on $R$, $L_{R}%
$, $M_{0}$, $L_{0}$,
$T$, $\alpha$, $\beta$, $\Lambda_{\alpha}\left(
B^{H}\right)  $.

Since $|x|^{2}=1$, then from the equation (\ref{equ4.1}) we clearly
have
\begin{equation}\label{equ4.2}
\left|\langle
2x,b(t,x)\rangle+\left[\int_{t}^{s}U(r)dr+\int_{t}^{s}V(r)dB^{H}(r)\right]\frac{1}{s-t}\right|=
\left|\langle
2x,\sigma(t,x)\rangle\frac{B^{H}(s)-B^{H}(t)}{s-t}\right|
\end{equation}
By (\ref{equ2}), there exists $\Omega_{0}\subset\Omega$ with
$\mathbb{P}(\Omega_{0})=1$ such that $\forall \omega\in\Omega_{0}$,
(\ref{equ4.1}) is satisfied and
\[\limsup_{t\rightarrow t_{0},~t\ge
t_{0}}\Big|\frac{B^{H}_{t}(\omega)-B^{H}_{t_{0}}(\omega)}{t-t_{0}}\Big|=+\infty.
\]
Let $\omega_{0}\in\Omega_{0}$. Then there is a subsequence
$r_{n}=r_{n}(\omega_{0})\downarrow t$ when $n\rightarrow\infty$,
such that
\begin{equation}\label{equ4.10}
\lim_{n\rightarrow\infty}\Big|
\frac{B^{H}_{r_{n}}(\omega_{0})-B^{H}_{t}(\omega_{0})}{r_{n}-t}\Big|=+\infty.
\end{equation}
Setting in \ref{equ4.2}
$s=r_{n}\wedge(t+\bar{h}(\omega_{0}))\in[t,t+\bar{h}(\omega_{0})]$
and passing to limit as $n\rightarrow\infty$, the left member,  via
Lemma \ref{lem4.1}, has limit $2\langle x,b(t,x)\rangle$.
Consequently, noting (\ref{equ4.10}), we must have
\[\langle x,\sigma(t,x)\rangle=0.
\]
and therefore
\[\langle x,b(t,x)\rangle=0.
\]
The proof of Lemma \ref{lem3.3} is complete.
\begin{flushright}
$\Box$
\end{flushright}
\textbf{Proof of Lemma \ref{lem3.4}}

Like the analysis in the proof of Lemma \ref{lem3.3}, the proof of Lemma \ref{lem3.4} is reduced to the following equivalent:

$\forall x$ such that $|x|=R$
\[(\langle 2x,b(t,x)\rangle,2x^{\ast}\sigma(t,x))\in S_{\varphi(K)}(t,|x|^2)\Leftrightarrow \langle x,b(t,x)\rangle\leq 0,
\quad \langle x,\sigma(t,x)\rangle=0
\]

$\forall x$, such that $|x|=r$,
\[(\langle 2x,b(t,x)\rangle,2x^{\ast}\sigma(t,x))\in S_{\varphi(K)}(t,|x|^2)\Leftrightarrow \langle x,b(t,x)\rangle\ge 0,
\quad \langle x,\sigma(t,x)\rangle=0.
\]

We only prove in the case $|x|=R$, the other one is similar.

\textit{Sufficient.} If $\langle x,b(t,x)\rangle\leq0, ~\langle x,\sigma(t,x)\rangle=0$,  taking $U(r)\equiv 0, ~V(r)\equiv 0$,
and we can choose $\bar{h}$ small enough such that
$\forall s\in[t,t+\bar{h}]$,
\[r^{2}\leq |x|^2+\int_{t}^{s}(\langle 2x,b(t,x)\rangle+U(r))dr+\int_{t}^{s}(2x^{\ast}\sigma(t,x)+V(r))dB^{H}(r)\leq R^{2},
\]
This means that $(\langle 2x,b(t,x)\rangle,\langle 2x,\sigma(t,x)\rangle)\in S_{\varphi(K)}(t,|x|^2)$.

\textit{Necessary.}  Since $(\langle 2x,b(t,x)\rangle,\langle
2x,\sigma(t,x)\rangle)\in S_{\varphi(K)}(t,|x|^2)$. Then there exist
random variable $\bar{h}=\bar{h}^{t,x}>0,$ and two stochastic
process
\[
\begin{array}[c]{ll}
U=U^{t,x}:& \Omega\times\left[  t,t+\bar{h}\right]  \rightarrow\mathbb{R} ,~\ U(t)=0
\\
V=V^{t,x}:& \Omega\times\left[  t,t+\bar{h}\right]  \rightarrow\mathbb{R} ,~\ V(t)=0
\end{array}
\]
such that for all $s,\tau\in\left[  t,t+\bar{h}\right]  $ and for
every $R>0$, $\left\vert x\right\vert \leq R:$
\[
\left\vert U\left(  s\right)  -U\left(  \tau\right)  \right\vert \leq
D_{R}\left\vert s-\tau\right\vert ^{1-\alpha},\quad\left\vert
V\left(  s\right)  -V\left(  \tau\right)  \right\vert \leq\tilde{D}_{R}\left\vert s-\tau\right\vert
^{min\{\beta,1-\alpha\}}%
\]
and
\[r^{2}\leq |x|^2+\int_{t}^{s}(\langle 2x,b(t,x)\rangle+U(r))dr+\int_{t}^{s}(\langle 2x,\sigma(t,x)\rangle+V(r))dB^{H}(r)\leq R^{2},
\]
where  $D_{R}$, $\tilde{D}_{R}%
$,  depend only on$R$, $L_{R}%
$, $M_{0}$, $L_{0}$,
$T$, $\alpha$, $\beta$, $\Lambda_{\alpha}\left(
B^{H}\right)  $.

Since $|x|=R$, then we yield
\begin{equation}\label{equ4.3}
\langle 2x,b(t,x)\rangle(s-t)+\langle 2x,\sigma(t,x)\rangle(B^{H}(s)-B^{H}(t))+\int_{t}^{s}U(r)dr+\int_{t}^{s}V(r)dB^{H}(r)\leq0
\end{equation}
By (\ref{equ3}),  there exists $\Omega_{0}\subset\Omega$ with
$\mathbb{P}(\Omega_{0})=\frac{1}{2}$ such that for each
$\omega_{0}\in\Omega_{0}$, (\ref{equ4.3}) is satisfied and there is
a sequence $t\leq s_{n}=s_{n}(\omega_{0})\leq
t+\bar{h}(\omega_{0})$, $s_{n}\downarrow t$, such that
\begin{equation}\label{equ4.11}
\lim_{s_{n}\downarrow
t}\frac{B^{H}_{s_{n}}(\omega_{0})-B^{H}_{t}(\omega_{0})}{s_{n}-t}=+\infty.
\end{equation}
Then we have
\begin{equation}\label{equ4.4}
\langle2x,\sigma(t,x)\rangle\frac{B^{H}_{\omega_{0}}(s_{n})-B^{H}_{\omega_{0}}(t)}{s_{n}-t}+
\left[\int_{t}^{s_{n}}U(r)dr+\int_{t}^{s_{n}}V(r)dB^{H}_{\omega_{0}}(r)\right]\frac{1}{s_{n}-t}
\leq -\langle 2x,b(t,x)\rangle.
\end{equation}
By Lemma \ref{lem4.1}
\[
\left[\int_{t}^{s_{n}}U(r)dr+\int_{t}^{s_{n}}V(r)dB^{H}_{\omega_{0}}(r)\right]\frac{1}{s_{n}-t}\rightarrow
0.
\]
and noting (\ref{equ4.11}), we derive that
\[\langle x,\sigma(t,x)\rangle\leq0.
\]
Similarly we  can prove $\langle x,\sigma(t,x)\rangle\ge0$, choosing
$\omega_{0}^{\prime}$ and a sequence $t\leq
r_{n}=r_{n}(\omega_{0}^{\prime})\leq t+\bar{h}(\omega_{0}^{\prime})$
and $r_{n}\downarrow t$ such that $\lim_{r_{n}\downarrow t}
\frac{B^{H}_{r_{n}}(\omega_{0}^{\prime})-B^{H}_{t}(\omega_{0}^{\prime})}{r_{n}-t}=-\infty$.
So
\[\langle x,\sigma(t,x)\rangle=0.
\]
Then from (\ref{equ4.3}), we have
\[\langle
2x,b(t,x)\rangle+\left[\int_{t}^{s}U(r)dr+\int_{t}^{s}V(r)dB^{H}(r)\right]\frac{1}{s-t}\leq0.
\]
and passing to limit $s\rightarrow t$, it follows, via Lemma
\ref{lem4.1},
\[\langle x,b(t,x)\rangle\leq0
\]
The proof of Lemma \ref{lem3.4} is finished.
\begin{flushright}
$\Box$
\end{flushright}
\textbf{Proof of Lemma \ref{lem3.5}} It is very similar to the proof
of Lemma \ref{lem3.4}, therefore we omit it.
\begin{flushright}
$\Box$
\end{flushright}

\section{Proof of the Comparison Theorem}\label{sec5}
\textbf{Proof of Theorem \ref{th3.10}}

We write  the two dimensional decoupled system%
\[
\left\{
\begin{array}
[c]{l}%
X_{s}^{t,x}=x+{\displaystyle\int_{t}^{s}}(b_{1}(r,X_{r}^{t,x}))dr+{\displaystyle\int_{t}^{s}}(\sigma_{1}(r,X_{r}^{t,x}))dB^{H}_{r},
\smallskip\;s\in[t,T] \medskip\\
Y_{s}^{t,y}=y+{\displaystyle\int_{t}^{s}}(b_{2}(r,Y_{r}^{t,y}))dr+{\displaystyle\int_{t}^{s}}(\sigma_{2}(r,Y_{r}^{t,y}))dB^{H}_{r},
\smallskip\;s\in[t,T]
\end{array}
\right.
\]
as
\[Z_{s}^{t,z}=z+{\displaystyle\int_{t}^{s}}(b(r,Z_{r}^{t,z}))dr+{\displaystyle\int_{t}^{s}}(\sigma(r,Z_{r}^{t,z}))dB^{H}_{r},\smallskip\;s\in[t,T]
\]
where
\[Z_{s}^{t,z}={X_{s}^{t,x}\choose Y_{s}^{t,y}}, ~z={x\choose y},~ b(r,Z_{r}^{t,z})={b_{1}(r,X_{r}^{t,x})\choose b_{2}(r,Y_{r}^{t,y})},
~\sigma(r,Z_{r}^{t,z})={\sigma_{1}(r,X_{r}^{t,x})\choose \sigma_{2}(r,Y_{r}^{t,y})}.
\]
we take $\varphi(z)=\varphi(x,y)=y-x$, then for every
$z\in\mathbb{R}^{2}$, $\varphi^{\prime}(z)=(-1,1)$, and
\[
\varphi^{\prime}(z)^{+}=\frac{1}{2}{-1 \choose1}
\]
so $\varphi\in\mathcal{H}$ and if we set $K=\{(x,y)\mid y-x\ge0\}$,
we have $\varphi(K)=\mathbb{R}^{+}$ and
$\varphi^{-1}(\varphi(K))=K$, then using the same method of the
proof of Theorem \ref{th3.2} we have
\[(b(t,z),\sigma(t,z))\in S_{K}(t,x)\Leftrightarrow (\varphi^{\prime}(z)b(t,z),\varphi^{\prime}(z)\sigma(t,z))\in S_{\mathbb{R}^{+}}(t,\varphi(z)).
\]
In fact, considering Theorem \ref{th3.1}, we only need prove
\[(\varphi^{\prime}(z)b(t,z),\varphi^{\prime}(z)\sigma(t,z))\in S_{\mathbb{R}^{+}}(t,\varphi(z))\Rightarrow (b(t,z),\sigma(t,z))\in S_{K}(t,x).
\]
Since
\[(\varphi^{\prime}(z)b(t,z),\varphi^{\prime}(z)\sigma(t,z))\in S_{\mathbb{R}^{+}}(t,\varphi(z))
\]
then for $z\in K$, there exist  a random variable
$\bar{h}=\bar{h}^{t,z}>0,$ and two stochastic process,
\[
\begin{array}
[c]{ll}
U_{1}=U_{1}^{t,z}:& \Omega\times\left[  t,t+\bar{h}\right]  \rightarrow\mathbb{R} ,~\ U_{1}(t)=0
\\
V_{1}=V_{1}^{t,z}:& \Omega\times\left[  t,t+\bar{h}\right]  \rightarrow\mathbb{R} ,~\ V_{1}(t)=0
\end{array}
\]
such that for all $s,\tau\in\left[  t,t+\bar{h}\right]  $ and for every $R>0$ and $\left\vert z\right\vert \leq R,$

\[
\left\vert U_{1}\left(  s\right)  -U_{1}\left(  \tau\right)  \right\vert \leq
D_{R}\left\vert s-\tau\right\vert ^{1-\alpha},\quad\left\vert
V_{1}\left(  s\right)  -V_{1}\left(  \tau\right)  \right\vert \leq\tilde{D}_{R}\left\vert s-\tau\right\vert
^{min\{\beta,1-\alpha\}}%
\]
and
\[
\varphi(z)+\int_{t}^{s}(\varphi^{\prime}(z)b(t,z)+U_{1}(r))dr+ \int_{t}^{s} (\varphi^{\prime}(z)\sigma(t,z)+V_{1}(r))dB^{H}(r) \in \varphi(K(s)).
\]
Let
\begin{eqnarray*}
f(r,y) & = &  \varphi^{\prime}(y)^{+}\left[U_{1}(r)-(\varphi^{\prime}(y)-\varphi^{\prime}(z))b(t,z)\right]\\
g(r,y) & = &
\varphi^{\prime}(y)^{+}\left[V_{1}(r)-(\varphi^{\prime}(y)-\varphi^{\prime}(z))\sigma(t,z)\right].
\end{eqnarray*}
It's obviously  that $f(r,y)$, $g(r,y)$ are independent of $y$.
%satisfy the conditions such that that the following SDE
Let
\[
\xi_{s}=z+\int_{t}^{s}(b(t,z)+f(r,\xi_{r}))dr+\int_{t}^{s}(\sigma(t,z)+g(r,\xi_{r}))dB^{H}(r), ~s\in[t,t+\bar{h}], ~|z|\leq R.
\]
Then we take
\begin{eqnarray*}
U(r) & = &  \varphi^{\prime}(\xi_{r})^{+}\left[U_{1}(r)-(\varphi^{\prime}(\xi_{r})-\varphi^{\prime}(z))b(t,z)\right]\\
V(r) & = & \varphi^{\prime}(\xi_{r})^{+}\left[V_{1}(r)-(\varphi^{\prime}(\xi_{r})-\varphi^{\prime}(z))\sigma(t,z)\right]
\end{eqnarray*}
According to the fractional It\^{o} formula, we have for all
$s\in\left[ t,t+\bar{h}\right]  $
\[
\varphi\Big(z+\int_{t}^{s}(b(t,z)+U(r))dr+\int_{t}^{s}(\sigma(t,z)+V(r))dB^{H}(r)\Big)
\]
\[=
\varphi(z)+\int_{t}^{s}(\varphi^{\prime}(z)b(t,z)+U_{1}(r))dr+\int_{t}^{s}(\varphi^{\prime}(z)\sigma(t,z)+V_{1}(r))dB^{H}(r)\in\mathbb{R}^{+}
\]
Clearly that
\[U(t)=0,\quad V(t)=0.
\]
Since for every $z\in\mathbb{R}^{2}$, $\varphi^{\prime}(z)=(-1,1)$,
$\varphi^{\prime}(z)^{+}=\frac{1}{2}{-1 \choose1}$ and together with
($\mathbf{H}_{1}$), ($\mathbf{H}_{2}$), it is clear that
\[
\left\vert U\left(  s\right)  -U\left(  \tau\right)  \right\vert \leq
\theta\left\vert s-\tau\right\vert ^{1-\alpha},\quad\left\vert
V\left(  s\right)  -V\left(  \tau\right)  \right\vert \leq\tilde{\theta}\left\vert s-\tau\right\vert
^{min\{\beta,1-\alpha\}}.%
\]
This means that
\[(b(t,z),\sigma(t,z))\in
S_{\varphi^{-1}(\varphi(K))}(t,x)=S_{K}(t,x).
\]

Just as Corollary \ref{cor2.7} said, if we want to get the
conditions for the viability of $K$, we only need to think about the
starting point $x\in\partial K$. Then the comparison theorem is
equivalent to prove that for all $t\in[0,T]$ and for any
$z={x\choose y}$ such that $x=y$, and $|z|\leq R$,
\[(\varphi^{\prime}(z)b(t,z),\varphi^{\prime}(z)\sigma(t,z))\in S_{\mathbb{R}^{+}}(t,\varphi(z))\Leftrightarrow b_{1}(t,x)\leq b_{2}(t,y),
\sigma_{1}(t,x)=\sigma_{2}(t,y).
\]

\textit{Sufficient}.  If $b_{1}(t,x)\leq b_{2}(t,y), \sigma_{1}(t,x)=\sigma_{2}(t,y)$, for $x=y$, we can take $U(r)\equiv 0, V(r)\equiv 0$,
and we have $\forall s\in[t,t+\bar{h}]$, and $z={x\choose y}$, such that $x=y$,
\[y-x+\int_{t}^{s}(\varphi^{\prime}(z)b(t,z)+U(r))dr+\int_{t}^{s}(\varphi^{\prime}(z)\sigma(t,z)+V(r))dB^{H}(r)=(b_{2}(t,y)- b_{1}(t,y))(s-t)\ge0,
\]
This means that $(\varphi^{\prime}(z)b(t,z),\varphi^{\prime}(z)\sigma(t,z))\in S_{\mathbb{R}^{+}}(t,\varphi(z))$.

\textit{Necessary}.  Since $(\varphi^{\prime}(z)b(t,z),\varphi^{\prime}(z)\sigma(t,z))\in S_{\mathbb{R}^{+}}(t,\varphi(z))$,
then there exist random variable $\bar{h}=\bar{h}^{t,z}>0,$ and two stochastic process
\[
\begin{array}[c]{ll}
U=U^{t,z}:& \Omega\times\left[  t,t+\bar{h}\right]  \rightarrow\mathbb{R} ,~\ U(t)=0
\\
V=V^{t,z}:& \Omega\times\left[  t,t+\bar{h}\right]  \rightarrow\mathbb{R} ,~\ V(t)=0
\end{array}
\]
such that for all $s,\tau\in\left[  t,t+\bar{h}\right]  $ and for every $R>0$ and
$\left\vert z\right\vert \leq R:$
\[
\left\vert U\left(  s\right)  -U\left(  \tau\right)  \right\vert \leq
D_{R}\left\vert s-\tau\right\vert ^{1-\alpha},\quad\left\vert
V\left(  s\right)  -V\left(  \tau\right)  \right\vert \leq\tilde{D}_{R}\left\vert s-\tau\right\vert
^{min\{\beta,1-\alpha\}}%
\]
and
\[y-x+\int_{t}^{s}((b_{2}(t,y)-b_{1}(t,x))+U(r))dr+\int_{t}^{s}((\sigma_{2}(t,y)-\sigma_{1}(t,x))+V(r))dB^{H}(r)\ge0,
\]
Since $y=x$, then we get
\begin{equation*}
(b_{2}(t,x)-b_{1}(t,x))(s-t)+(\sigma_{2}(t,x)-\sigma_{1}(t,x))(B^{H}(s)-B^{H}(t))+\int_{t}^{s}U(r)dr+\int_{t}^{s}V(r)dB^{H}(r)\ge0.
\end{equation*}
With the same analysis in the proof of Lemma \ref{lem3.4}, we obtain
that for every $R>0$
\[b_{1}(t,x)\leq b_{2}(t,x),\quad \sigma_{1}(t,x)=\sigma_{2}(t,x), \quad\forall |x|\leq
R.
\]
This complete the proof of Comparison Theorem.
\begin{flushright}
$\Box$
\end{flushright}
\textbf{Acknowledgement}\quad The authors express special thanks
to Rainer Buckdahn and Lucian Maticiuc for their useful
suggestions and discussions.

\end{document}